%Paper: q-alg/9509018
%From: gannon@mpim-bonn.mpg.de (Terry Gannon)
%Date: Sun, 17 Sep 1995 16:33:49 +0200
%Date (revised): Tue, 26 Sep 1995 09:46:08 +0100

\magnification=\magstep1
\font\huge=cmr10 scaled \magstep2    \font\smal=cmr7
\font\smit=cmmi7
\overfullrule=0pt

%%%%% Blackboard bold characters %%%%%
  \catcode`\@=11
  \font\tenmsa=msam10 \font\sevenmsa=msam7 \font\fivemsa=msam5
  \font\tenmsb=msbm10 \font\sevenmsb=msbm7 \font\fivemsb=msbm5
  \newfam\msafam        \newfam\msbfam
  \textfont\msafam=\tenmsa  \scriptfont\msafam=\sevenmsa
  \textfont\msbfam=\tenmsb  \scriptfont\msbfam=\sevenmsb
  \scriptscriptfont\msafam=\fivemsa  \scriptscriptfont\msbfam=\fivemsb
  \def\hexnumber@#1{\ifcase#1 0\or1\or2\or3\or4\or5\or6\or7\or8\or9\or
        A\or B\or C\or D\or E\or F\fi }
   \def\Bbb{\ifmmode\let\next\Bbb@\else
  \def\next{\errmessage{Use \string\Bbb\space only in math mode}}\fi\next}
  \def\Bbb@#1{{\Bbb@@{#1}}} \def\Bbb@@#1{\fam\msbfam#1}
%%%%%         (end of Blackboard bold characters)

\def\sp{\;}    \def\eg{{\it e.g.}$\sp$} \def\ie{{\it i.e.}$\sp$}
  \def\Q{{\Bbb Q}}   \def\L{{\cal L}}
\def\M{{\cal M}}             \def\Z{{\cal Z}} \def\om{\omega}
\def\g{\hat{g}}  \def\la{\lambda} \def\si{\sigma} \def\eps{\epsilon}
\def\z{{\Bbb Z}} \def\GAL{{\rm Gal}({\Q}^{ab}/\Q)}  \def\ga{\gamma}

\def\Bi{[1]}  \def\Ro{[2]}   \def\Jo{[3]} \def\BL{[4]}
\def\W{[5]} \def\RT{[6]} \def\Kac{[7]}     \def\CG{[8]}
\def\AC{[9]} \def\Ga{[10]}  \def\FSS{[11]} \def\GJW{[12]} \def\ST{[13]}
\def\BC{[14]} \def\Ver{[15]}    \def\BN{[16]}      \def\V{[17]}
 \def\Wi{[18]}  \def\Ma{[19]} \def\dB{[20]} \def\RTW{[21]} \def\Gu{[22]}
\def\GH{[23]}  \def\Gr{[24]} \def\Dr{[25]}  \def\De{[26]} \def\Ih{[27]}
\def\BR{[1,2]}  \def\BMW{[18,19]}

{\nopagenumbers
\rightline{August, 1995}\bigskip\bigskip
\centerline{{\bf \huge  Galois Relations on Knot Invariants}}
\bigskip \bigskip

\centerline{Terry Gannon\footnote{${}^{\dag}$}{{\smal Address from
Sept 1995-Aug 1996: Max-Planck-Institut f\"ur Math,
Bonn, D-53225}} \footnote{}{{\smal Permanent address as of Sept 1996:
Math Dept, York Univ, North York, M3J 1P3, Canada}} and Mark A.\
Walton\footnote{${}^{\ddag}$}{{\smal Physics Dept, University of Lethbridge,
Lethbridge, Alberta, T1K 3M4,
Canada}}}
\bigskip\bigskip

\centerline{{\bf Abstract}}
\medskip

\noindent We discuss the existence of  Galois relations obeyed by certain
link invariants. Some of these relations have recently been identified and
exploited
within the context of conformal field theory and Lie/Kac-Moody representation
theory. These relations should aid in computing knot invariants.
They probably have an interpretation in terms of quasitriangular (quasi) Hopf
algebras. They could also have a topological interpretation, and may serve as
a concrete model for related ideas of Degiovanni, Drinfeld, Grothendieck,
Ihara, and others.
\vfill\eject}\pageno=1

\noindent{\bf 1. Introduction}\medskip

The classification and study of knots is an old
problem \BR. One of the most fruitful ideas  has been to assign to each knot
a set of objects (e.g.\ numbers or polynomials): two knots will be equivalent
only if these objects are identical. The first  knot polynomial is due to
Alexander in 1928, but
many more invariants have been discovered since (see \eg \Jo). In fact there
is now
an embarrassment of riches, but it is still an open question whether all the
known invariants will completely distinguish all knots.

Three important goals in the subject are: to find convenient aids for
computing these invariants; to obtain a topological understanding
of these invariants; and to bring some order to the multitude of invariants
by finding systematic relationships among them. The traditional approach to
computation has been to use some sort of skein relation, but these are
practical only for the simplest invariants and knots. An adequate topological
understanding is still missing \Bi\ for most invariants. Most invariants are
redundant; the most significant
relationship known among them is the theorem of Birman-Lin
\BL\ that the Vassiliev invariants contain the Lie algebra ones.

In the program introduced in this paper, we hope
to provide new means to approach the first two goals, by using the Galois
group $\GAL$ (by $\Q^{ab}$ we mean a maximal abelian extension of $\Q$)
to directly address the latter.
In \W\ Witten establishes connections between 3-dimensional
topology and 2-dimensional rational conformal field theory (RCFT); we will
exploit that connection as we generalize relations recently obtained in RCFT.
For concreteness we will focus here on the Witten link invariants \W,
which equal the quantum group link invariants \RT\ for
the quantum enveloping algebras $U_q(g)$, where $g$ is a Lie algebra and $q$
a root of unity.
The special property possessed by these invariants is the presence of the
 Lie algebra -- it makes the appearance of Galois relations more explicit.
Later in this paper we will discuss the possibility that
there is an action of the group $\GAL$ on knots themselves,
and its possible relation to certain ideas of Grothendieck and others.

Each (directed) link  $\L$ in each closed, connected and orientable
3-manifold\footnote{$^1$}{{\smal
To be well-defined, the Witten invariant requires the 3-manifold to be
{\smit framed} \W, i.e.\ a trivialization of its tangent bundle must be
chosen. A different choice however merely multiplies the invariant by a certain
root of unity. In this paper we can and will safely ignore this complication.}}
$\M$ (\ie an
embedding $(S^1)^t\rightarrow \M$) possesses several Witten invariants.
Namely, choose any nontwisted affine Kac-Moody algebra $\g=X_r^{(1)}$, any
integer $k>0$, and any level $k$ highest weights
$\la^{1},\ldots,\la^{t}\in P_+^k(\g)$  (we will review  affine
algebra representation theory in the following section; see also \Kac). Then
there is a Witten invariant  \W\
$$\Z(\L,\g,\{\la^{i}\})\in{\Bbb C}\ .\eqno(1)$$
We will discuss eq.(1) in section 2; these numbers $\Z$
(topological invariants for $\L$) are the starting point for this paper.

We will find a Galois relation mixing some of these $\Z$, which in some cases
will be non-linear, and can involve different
$\L$ (and even $\M$). This project is motivated in part by \CG, in which was
found a simple Galois relation\footnote{$^2$}{{\smal
More generally, this holds for any RCFT \CG,
which suggests that the comments in this paper should also hold in non-Lie
theoretic contexts -- for example the quasi-Hopf link invariants \AC\ --
though perhaps at the expense of ease of computation.}}
for the Kac-Peterson matrix $S$ and corresponding
fusion coefficients $N_{\la,\mu}^\nu$ (see eqs.(3)).
$S_{\la,\mu}$ and $N_{\la,\mu}^\nu$ turn out to be examples of link invariants,
so it is natural to ask if other link invariants will also obey reasonably
simple Galois relations. This is made much more plausible by the observation
that the Witten invariants can be obtained from surgery (which involves $S$)
and skein relations (which involve \eg roots of unity).
This is reviewed and generalized in section 3, where we also make some remarks
about what the general Galois relations will look like and how they can be
found.

The Galois symmetry on $S$ has been very effective in the classification
of RCFTs (see \eg \Ga). The action on $N$ has been used \FSS\ to construct
new exceptional partition functions of RCFTs. They were used in \GJW\ to
obtain relations among weight multiplicities for representations of
finite-dimensional Lie algebras. Further applications\footnote{$^3$}{{\smal
Galois has made other appearances in RCFT -- see e.g.\ \ST.}} are discussed
in \BC.
Thus it is certainly  possible that the generalization of these relations to
other link invariants will yield other valuable formulae, with applications
lying outside the  context of RCFT. This should
 be particularly true once a more topological understanding of
these invariants is found. In the final section we address some of these
points.

A well-known collection of relationships among the Witten invariants
involves the transformation
$\la\mapsto C\la$, taking a Lie algebra highest weight to the one
contragredient to it.
Witten's invariants are preserved if $C$ is applied simultaneously to all
weights $\la^{i}$:
$$\Z(\L,\g,\{C\la^{i}\})=\Z(\L,\g,\{\la^{i}\})\ .\eqno(2a)$$
This has a topological meaning: the Witten invariants of a link $\L$ are
equal to those of its {\it reverse} $\L^r$:
$$\Z(\L,\g,\{\la^{i}\})=\Z(\L^r,\g,\{\la^{i}\})\ .\eqno(2b)$$
$\L^r$ is obtained from $\L$ by switching the orientations of all of
its components. Since some knots are inequivalent to their reverse,
eq.(2b) implies that the Witten invariants cannot distinguish all
inequivalent knots.

$C$ also has a privileged role in RCFT, but in \CG\
it is found that $C$ is only one among a whole family of Galois
``symmetries'' $\si\in\GAL$ which obey similar equations
($C$ corresponds to complex conjugation); whenever $C$ possesses some
property in RCFT, an analogue is usually possessed by all other $\si$. In
this paper we propose that a similar situation should hold in knot theory,
\ie that analogues of eq.(2a) and possibly eq.(2b) hold for all $\si$.

In this short note we can do little more than briefly state the general ideas,
give some simple examples, and speculate about possible future
developments. Much more work is required before their scope and
usefulness can be uncovered. However,
we believe the possibilities are intriguing enough to justify
publication of this work at this early stage.

\vfill\eject
\bigskip\noindent{\bf 2. Witten link invariants}\medskip

We begin by briefly describing some objects in affine algebra representation
theory (see \Kac\ for more details),
then proceed to give a few simple facts about links (see also \BR), and end by
reviewing the basic theory of Witten invariants (see \eg \W).

There is a one-to-one correspondence between finite-dimensional semi-simple Lie
algebras $X_r$, and nontwisted affine Kac-Moody algebras $\g=X_r^{(1)}$.
The integrable highest weight representations of $\g$ are characterized
by their highest weights $\mu=\sum_{i=0}^r\mu_iw^i$, where each $\mu_i
\in \z_{\ge}$ (the $w^i$ are the fundamental weights). The set of all highest
weights for $\g$ are partitioned into finite subsets $P_+^k(\g)$, indexed by
the {\it level} $k$ (\eg for $\g=A_r^{(1)}$, $k=\sum_{i=0}^r\mu_i$). The
characters $\chi_\mu$ of $\g$, for $\mu\in P_+^k(\g)$, define a
representation of the modular group SL$_2(\z)$:
$$\eqalignno{\chi_\mu(\tau+1,z,u)=&\sum_{\nu\in P_+^k} T_{\mu,\nu}\,\chi_\nu
(\tau,z,u)&(3a)\cr
\chi_\mu({-1\over \tau},{z\over \tau},u-{(z\,|\,z)\over \tau})=&
\sum_{\nu\in P_+^k} S_{\mu,\nu}\,\chi_\nu(\tau,z,u)\ .& (3b)\cr}$$
$S$ and $T$ are called the {\it Kac-Peterson matrices}. They are both unitary.
$T$ is diagonal, while the elements of $S$ are related to certain values of
the characters of $X_r$. The map $C$ appearing in eq.(2a) corresponds to the
permutation matrix $S^2$.

We are also interested in the {\it fusion coefficients} $N_{\la,\mu}^\nu$,
which  are defined by Verlinde's formula \Ver:
$$N_{\la,\mu}^\nu=\sum_{\ga\in P_+^k}S_{\la,\ga}\,{S_{\mu,\ga}\over
S_{k\om^0,\ga}}\,S^*_{\nu,\ga}\ .\eqno(3c)$$
The $N_{\la,\mu}^\nu$ will always be non-negative integers, and are
closely related to the tensor product coefficients of $X_r$. They have
appeared in several different contexts in recent years, including of course
RCFT.

A {\it link} is an embedding of $S^1\times\cdots\times S^1$
into a 3-manifold $\M$.
To avoid certain wild knots (\eg infinite connected sums),
additional conditions (\eg smoothness or polygonality) are required \BR,
but will not be given here. All links here are directed.
A {\it knot} is a link with precisely one connected component.

A {\it link invariant} is a function  which assigns to each
link $\L$ an object $f(\L)$ in such a way that equivalent links are
assigned equivalent (usually equal, or equal up to some factor) objects.
A simple example is the number of connected components in the link.

We call two links $\L,\L'$ equivalent when they are {\it ambient isotopic}
\Ro, \ie  if there exists an
orientation-preserving homeomorphism $h:\M\rightarrow \M$ with $h(\L)=\L'$.
In $S^3$ this is the same as saying  their link
diagrams (\ie projections) can be related by some sequence of the
familiar Reidemeister moves I, II, III \Bi. Two links are called {\it regular
isotopic} when their diagrams are invariant under Reidemeister II, III.
In the literature both ambient and regular isotopic link
invariants are discussed: the Jones polynomial and its generalization
HOMFLY \Jo, as well as Conway's normalization of the Alexander polynomial,
are ambient; while
ribbon invariants, Kauffman's polynomial, and Witten's invariants are
regular. Reidemeister I says that a localized loop in a string can be
straightened; the special thing about this is that it changes the
``self-winding number'' of the knot. In the case of the Kauffman polynomial
and  the Witten invariants, ambient isotopy can always be recovered by
introducing a phase factor depending on this self-winding number.
The problem is that the resulting ambient isotopic invariant will depend
on an arbitrary choice of overall phase factor, reflecting the fact that
different ribbons can have the same knot as their ``spine''.

In the language of quantum field theory, the Witten invariant in eq.(1) is the
unnormalized expectation value of Wilson lines $\L$ in the Chern-Simons
theory whose gauge group is the compact exponentiation of $X_r$. The
Lagrangian is proportional  to $k$, and the components of $\L$ carry the
representations of highest weights $\overline{\la^i}=(\la_1^i,\ldots,\la_r^i)$.
Mathematically, $\Z$ is definable using the
Reshetikhin-Turaev functor $F$ \RT, which relates a category of directed
ribbon graphs
coloured by weights in $P_+^k(\g)$ to a category of finite-dimensional
representations of $U_q(X_r)$. An alternative is the approach
given in \BN\ for Lie algebra invariants, using ``Chinese character diagrams''.

Some (but not all) of the numbers in eq.(1) for fixed $\L$ and $\g$ may be
collected and reinterpreted as polynomials in a formal parameter $q$,
evaluated at the root of unity $q=\exp[2\pi i/(k+h^\vee)]$, where $h^\vee$
is  the {\it dual Coxeter number} of $X_r$.
For example, Jones' original polynomial $J(z)$, originally defined only
for $\L\subset S^3$, is related to
$\Z(\L,{A_1}^{(1)},w^1+(k-1)w^0)$, where every
strand of $\L$ is coloured by the (level $k$ ``lift'' of the) fundamental
representation $w^1$ of $A_1$. In particular,
this Witten invariant $\Z$  will equal $\omega\,J(\exp[2\pi i/(k+2)])$,
where $\omega$ is some appropriately chosen root of unity ($J$ is
uniquely determined by its value at those infinitely many
points). HOMFLY turns out to be the generalization of the Jones
polynomial to ${A_r}^{(1)}$, but again the only representation which appears
there  is the first fundamental one of $A_r$.

Witten invariants are a subset of the quantum group
invariants \RT, which turn out \BL\ to be special cases of the Vassiliev
invariants \V. It is not yet known however whether Vassiliev invariants
are actually more powerful than Witten invariants; both can be graded by a
positive integer called {\it degree}, and for degree up to 9 the dimensions of
the spaces spanned by
these invariants remain equal \BN. It is also unknown whether
Vassiliev invariants distinguish all knots.

One of the main problems is how to compute the $\Z$, even in principle.
Usually two steps are followed. The first is to replace $\M$ with a more
convenient
manifold, like $S^3$ or $S^1\times R^h$ ($R^h$ is the genus $h$ Riemann
surface). The main way to do this is surgery about an appropriately chosen
knot. The result \W\ is
an expression for the original $\Z$ as a linear combination of $\Z$'s
for a new link; the coefficients
of the linear combination are computable using the Kac-Peterson
matrices $S$ and $T$. In some special cases,
connected and disconnected sums of links are also effective.

So because of surgery we can compute any $\Z(\L,\g,\{\la^{i}\})$ if we know
all $\Z(\L'\subset S^3,\g,\{\mu^{j}\})$, say.
The problem then reduces to working in $\M=S^3$, and trying to simplify
the link $\L$ to some link we can handle. The main way to do this involves
skein relations, along with an expression saying what phase is picked up
when Reidemeister I is performed on the knot. These relations are given
in \BMW, and in general require a generalization \Wi\ of eq.(1)
from embeddings in $\M$ of $(S^1)^t$,
to embeddings in $\M$ of disjoint unions of {\it trivalent graphs} -- \ie
graphs with three
segments at each vertex. Each segment is assigned a weight, and $\Z$ will
vanish unless the fusion coefficient at each vertex is non-zero.   Ref.\
\BN\ gives an alternate but related approach.

The skein relations permit one to reduce the computation of the $\Z(\L)$
to the evaluation of invariants of ``tetrahedron
graphs'' (trivalent graphs with 4 vertices and 6 sides). The coefficients
of these skein relations are awkward: they involve the square-roots
$\sqrt{S_{0,\la}/S_{0,0}}$, along with ``braiding matrix elements''
(no explicit expression for these are known in general) and $T_{\la,\la}$. In
\Ma\ it is suggested that these tetrahedron invariants can be computed
as the solutions of linear equations (consistency conditions from the skein
relations). However there is no proof these equations can always be
inverted, and in any event this approach is very impractical except in the
simplest cases.

For use in the following section, we give here a few of the simpler
Witten invariants (for readability we will write `0' for the weight
`$k\om^0$').

\item{(i)} \quad $t$ parallel (\ie unlinked) unknots, in $S^3$ (independent
of orientation):
$$\Z=:{\cal D}_{\la^{1},\ldots,\la^{t}}=S_{0,0}\prod_{i=1}^t
{S_{0,\la^{i}}\over S_{0,0}}\ ;\eqno(4a)$$

\item{(ii)} \quad $t$ parallel unknots $S^1\times\{p_i\}$, $1\le i\le t$, in
$S^1\times R^h$ (oriented the same way): this is the Verlinde
dimension \Ver\
$$\Z={\cal V}^{h,t}_{\la^{1},\ldots,\la^{t}}
=\sum_{\mu\in P_+^k} (S_{0,\mu})^{2(1-g)}{S_{\la^{1},\mu}\over S_{0,\mu}}
\cdots{S_{\la^{t},\mu}\over S_{0,\mu}}\ ;\eqno(4b)$$

\item{(iii)} \quad a chain of $t$ linked unknots, in $S^3$ (oriented in the
same way):
$$\Z=:{\cal
C}_{\la^{1},\ldots,\la^{t}}=S_{0,\la^1}\prod_{i=1}^{t-1}{S_{\la^i,\la^{i+1}}
\over S_{0,\la^{i}}}\ ; \eqno(4c)$$

\item{(iv)} \quad a central unknot in $S^3$ with weight $\la^{0}$, with $t$
unknots linked around it (like keys around a key chain), all oriented the same:
$$\Z=:{\cal
S}_{\la^{0};\la^{1},\ldots,\la^{t}}=S_{0,\la^{0}}\prod_{i=1}^t
{S_{\la^{0},\la^{i}}\over S_{\la^{0} ,0}}\ .\eqno(4d)$$

\medskip\noindent Notice that ${\cal
S}_{0;\la^{1},\ldots,\la^{t}} = {\cal D}_{\la^{1},\ldots,\la^{t}}.$
 Also, $S_{\la,\mu}={\cal S}_{\la;\mu}={\cal C}_{\la,\mu}$ and $N_{\la,\mu}^\nu
={\cal V}_{\la,\mu,C\nu}^{0,3}$.

\medskip\bigskip\noindent{\bf 3. Galois relations for the Witten invariants}
\medskip

Next we review the presence of Galois in RCFT,
and conjecture how it will extend to the Witten link invariants.
What we are after is some kind of generalization of eq.(2a) to any Galois
automorphism $\si\in\GAL$ of the underlying RCFT.

The restriction to $\GAL$ is forced by the
invariants: they must lie in cyclotomic extensions of
$\Q$. For links in $S^3$, this can be seen directly from \BL, and since
it is true for the $S$ and $T$ matrices, by surgery it is true for links
in all other 3-manifolds $\M$.

It was discovered in \CG\ (based on work in \dB, and generalizing and
reinterpreting \RTW) that there is a Galois ``symmetry'' present in any
RCFT. In particular, let $M$ be the field extension of $\Q$ obtained
by adjoining to it all the matrix elements $S_{\la,\mu}$. It turns out that
$M$ will always lie in a cyclotomic field (for
RCFT based on an affine algebra, this follows directly from Kac-Peterson).
Choose any
$\si\in {\rm Gal}(M/\Q)$. Then there is a permutation $\la\mapsto \si\la$
of $P_+^k(\g)$, and a choice of signs $\eps_\si(\la)\in\{\pm 1\}$, such that
$$\si(S_{\la,\mu})=\eps_\si(\la)\,S_{\si\la,\mu}=\eps_\si(\mu)\,
S_{\la,\si\mu}\ .\eqno(5a)$$
Both $\eps_\si$ and $\la\mapsto \si\la$ have a geometric meaning in terms
of the affine Weyl group of $\g$ \CG, and can be readily calculated in
practice.
Eq.(5a) is the basic equation of \CG, and the source of the name ``Galois''
for the family of relations considered in this paper. From this and eq.(3c),
it is possible to derive \CG\ the expression
$$N^{\nu}_{\si\la,\si\mu}=\eps_\si(0)\,\eps_\si(\la)\,\eps_\si(\mu)
\sum_{\ga\in P_+^k}\eps_{\si}(\ga)\, N^{\ga}_{\la,\mu}\,
N^{\nu}_{\si\ga,\si 0}\ .\eqno(5b)$$

As we saw at the end of the last
section, both $S_{\la,\mu}$ and $N_{\la,\mu}^\nu$ can be interpreted as
link invariants, corresponding
respectively to the Hopf link (a chain of 2 linked unknots) in $S^3$, and
3 parallel unknots in $S^1\times S^2$. Since they
both satisfy a Galois relation, then perhaps many other invariants do too.
It is easy to find these relations when
explicit expressions exist for the invariants. For example, the same
calculation which gave us eq.(5b) can be used to give us:
$${\cal V}^{h,t}_{\si\la^{1},\ldots,\si\la^{t-1},\la^{t}}=
{\eps_\si(\la^1)\over \eps_\si(0)}\cdots{\eps_\si(\la^{t-1})\over\eps_\si(0)}
\sum_{\mu\in P_+^k}{\eps_\si(\mu)\over \eps_\si(0)}\,{\cal V}^{h,t}_{\la^{1},
\ldots, \la^{t-1},C\mu}\,{\cal V}^{0,2h+t}_{\si\mu,\si 0,\ldots,\si 0,
\la^{t}}\ .\eqno(5c)$$
{}From eqs.(4b) and (5a) we find that each
${\cal V}^{h,t}_{\la^1,\ldots,\la^t}$
is fixed by $\si$ -- in fact they are non-negative integers.
Galois relations for the  other invariants in eqs.(4) can also be
found. Simple generalizations of eq.(5a) are
$$\eqalignno{\si\left({\cal
S}_{\la^{0};\la^{1},\ldots,\la^{t}}\right) =& \,\epsilon_\si(\la^{0})\
{\cal S}_{\si\la^{0};\la^{1},\ldots,\la^{t}}\ &(5d)\cr
{\cal S}_{\si\la^0;\la^1,\ldots,\la^t}=&\prod_{i=0}^t{\eps_\si(\la^i)\over
\eps_\si (0)}\ {\cal S}_{\la^0;\si 0}\ {{\cal S}_{\la^0;\si\la^1,\ldots,
\si\la^t}\over {\cal S}_{\la^0;\si 0,\ldots,\si 0}}\ .&(5e)\cr}$$

By analogy with eqs.(5), we are looking for polynomial relations over $\z$
among link invariants.
The link invariants may not all correspond to the same link $\L$, or even the
same background manifold $\M$, but the weights appearing in them are either
the inputs $\la^i$, or $0$, or their images $\si\la^i$, $\si 0$ under some
$\si\in\GAL$ (as in eq.(5b), there can also appear dummy variables
$\ga$ involved in a sum -- a trace -- over all of $P_+^k$). In addition, in
some cases applying $\si$ directly to the invariant will yield a simple
expression.

There is a way, using the machinery developed in \GJW, to {\it linearize}
eq.(5b), using the dominant weight multiplicities $m_\la^\mu$ of
the highest weight module $L(\la)$ of $g$. In particular, we can write
eq.(5b) in the equivalent form
$$N_{\si\la,L_\si(\nu)}^{\si\mu}:={\|W\overline{\nu}\|\over
\|W\overline{\pi_\si(\nu)}\|}\, \sum_{\ga\in P_+^k}
\ell_{\pi_\si(\nu)}^{\,\ga}\,N_{\si\la,\ga}^{\si\mu}=
\eps_\si(\la)\,\eps_\si(\mu)\,
N_{\la,L_{id}(\nu)}^\mu\,
,\eqno(6a)$$
where $\pi_\si(\nu):=\si(\nu-\rho)+\rho,$ and the $\ell$'s are
the elements of the matrix $L=M^{-1}$ inverse to the matrix $M$ formed
 from the multiplicities $m_{\la}^\mu$ (both $L$ and $M$ are lower triangular,
with 1's along the diagonal).  $\|W\overline{\nu}\|$ is the order of the
(finite) Weyl orbit of $\overline{\nu}=(\nu_1,\ldots,\nu_r)$.
A similar trick works for other nonlinear Galois relations; for example,
$$
{\cal C}_{L_{id}(\la^{1}),\la^{2},\la^{3},L_{id}(\la^{4})}=
\eps_{\si^{-1}}(\la^{2})\,\eps_\si(\la^{3})
\,{\cal C}_{L_\si(\la^{1}),\si^{-1}\la^{2},\si\la^{3},L_{\si^{-1}}(\la^{4})}
\,\ .\eqno(6b)
$$
By the subscript $L_\si(\la)$ in eqs.(6) we mean the linear combination
over that position as defined on the r.h.s.\ of eq.(6a).

Eq.(6a) is linear in the fusion coefficients; the price is the introduction of
the ``inverse multiplicities'' $\ell_{\la}^\mu$. However these are easy to
compute (see \eg eqs.(2.6),(5.2) of \GJW), and so can be regarded
essentially
as known constants. This new relation (6a) should serve as a useful means of
computing fusions -- for example it directly provides an expression for the
fusion coefficients $N_{\si 0,\la}^\mu$ using the $\ell_\ga^\nu$'s and
the affine reflection $r_0$ built into $\pi_\si$ (this is where the
dependence of $N$ on $k$ enters).
Eq.(6a) together with the ``normalization'' $N_{0,0}^0=1$ may suffice in fact
to determine all other fusion coefficients.

Finding Galois relations for more complicated $\Z$'s will be difficult,
until we possess a more reliable way of computing these $\Z$'s.
If we use the skein relations, we presumably must know
Galois relations for \eg the tetrahedron graphs, perhaps by studying the
consistency conditions of \Ma. All this can probably be avoided in the
quantum group approach of \RT: in the most optimistic case we would be looking
in that language for a Galois action on the two categories (the representation
category, and the category of coloured directed ribbon graphs) which
commutes with the covariant functor $F$ between them. An interesting
alternative within the Chern-Simons theory comes from e.g.\ \Gu.

\bigskip\noindent{{\bf 4. Conclusion and speculations}}\medskip

Arguing by analogy with RCFT, we suggest that there could be a generalization
of eq.(2a) to the other automorphisms in $\GAL$ ($C$ in
eq.(2a) corresponds to complex conjugation). We provide some
simple examples supporting this claim. This
Galois ``symmetry'' has been useful in RCFT, so it may also be of
use in analyzing knot invariants. For example it should help in computing
some invariants. It would be interesting to find out if these Galois relations
exist merely because of the presence of Lie algebras in these Witten
invariants, or whether in some form they will persist for more general
link invariants (\eg Vassiliev). If they do persist, then perhaps they
will obey relations obtained from $\si\in{\rm Gal}(\overline{\Q}/\Q)$,
the absolute Galois group of $\Q$.

Hopf algebras and Galois theory have deep connections \GH. It is very
possible that the ultimate source of these Galois relations, both in
RCFT and for link invariants, lie in quasitriangular (quasi-)Hopf algebras.
This is then a natural direction to pursue the ideas in this paper, both
to clarify and to generalize them.

In eq.(2b) we find that complex conjugation has a topological interpretation.
It would be desirable to find one for the other elements of $\GAL$.
Unfortunately a more systematic treatment of these Galois relations, both
algebraically and topologically, will
be difficult without a more systematic algorithm for computing the Witten
invariants, or without a concrete understanding of the Galois actions
(if they exist) on the two Reshetikhin-Turaev categories.

There is a certain resemblance of these ideas to much more ambitious ones
sketched by Grothendieck \Gr, Drinfeld \Dr, Degiovanni \De, and others
(see in particular Ihara \Ih\ and references therein). Grothendieck proposed
to study ${\rm Gal}(\overline{\Q}/\Q)$  using the Teichm\"uller
tower formed from the moduli spaces $\M_{h,t}$ of Riemann surfaces of
genus $h$ with $t$ punctures; the elements of ${\rm Gal}(\overline{\Q}/\Q)$
act as outer automorphisms of the tower. Drinfeld related this
to a universal braid transformation group acting on structures of
quasitriangular quasi-Hopf algebras. An RCFT interpretation
has been suggested in \De: any RCFT may provide a projective  representation
of Grothendieck's tower, and there could be an action of ${\rm Gal}
(\overline{\Q}/\Q)$  on the data of RCFT (\eg the $S$ and $T$ matrices).
Perhaps the program introduced in this paper can serve as a concrete
realization of some aspects of these deep ideas.

\bigskip\noindent{{\bf Acknowledgements.}}\quad TG thanks Antoine Coste,
Pascal Degiovanni, and Ramin Naimi for several helpful conversations, and
the Math Dept at Concordia University for hospitality. MW acknowledges a
Research Grant from NSERC.

\bigskip\noindent{\bf References}\medskip

\item{1.}  Birman, J., {\it Bull.\ Amer.\ Math.\ Soc.}\ {\bf 28}, 253 (1993).

\item{2.} Rolfsen, D., {\it Knots and Links}, Publish or Perish Inc.,
Berkeley, 1976.

\item{3.} Jones, V., {\it Bull.\ Amer.\ Math.\ Soc.}\ {\bf 12}, 103 (1985);

\item{} Freyd, P.,  Yetter, D.,  Hoste, J.,  Lickorish, W.,  Millet, K.,
and Ocneanu, A., {\it Bull.\ Amer.\ Math.\ Soc.}\ {\bf 12}, 239 (1985).

\item{4.} Birman, J.\ S.\ and Lin, X.-S., {\it Invent.\ Math.}\  {\bf
111}, 225 (1993).

\item{5.} Witten, E., {\it Commun.\ Math.\ Phys.}\ {\bf 121}, 351 (1989).

\item{6.} Reshetikhin, N.\ and Turaev, V., {\it Commun.\ Math.\ Phys.}\
{\bf 127}, 1 (1990); {\it Invent.\ Math.}\ {\bf 103}, 547 (1991).

\item{7.} Kac, V.\ G., {\it Infinite Dimensional Lie Algebras},
Cambridge University Press, Cambridge, 1990.

\item{8.}  Coste, A.\ and Gannon, T., {\it Phys.\ Lett.}\ {\bf B323}, 316
(1994).

\item{9.} Altschuler, D.\ and Coste, A., {\it Commun.\ Math.\ Phys.}\
{\bf 150}, 83 (1992).

\item{10.} Gannon, T., {\it Commun.\ Math.\ Phys.}\ {\bf 161}, 233
(1994); ``The classification of SU(3) modular invariants revisited'',
{\it Annales de l'I.\ H.\ P.: Phys.\ Th\'eor.}\ (to appear).

\item{11.} Fuchs, J.,   Schellekens, A.\ N., and Schweigert, C.,
{\it Nucl.\ Phys.}\ {\bf B437}, 667 (1995).

\item{12.}  Gannon, T., Jakovljevic, C., and Walton, M.\ A.,
{\it J.\ Phys.}\ {\bf A28}, 2617 (1995).

\item{13.} Stanev, Ya.\ S.\ and Todorov, I.\ T., ``On Schwarz problem for
the $\widehat{SU}_2$ Knitzhnik-Zamolodchikov equation'', {\it Lett.\ Math.\
Phys.}\ (to appear).

\item{14.} Buffenoir, E., Coste, A.,  Lascoux, J.,  Buhot, A.,  and
Degiovanni, P., ``Precise study of some number fields and Galois actions
occurring in conformal field theory'', ENSLAPP-L-477/94.

\item{15.} Verlinde, E., {\it Nucl.\ Phys.}\ {\bf B300}, 360 (1988).

\item{16.} Bar-Natan, D., {\it Topology} {\bf 34}, 423 (1995).

\item{17.} Vassiliev, V.\ A., in {\it Theory of Singularities and its
Applications}, (V.\ I.\ Arnold, ed.), Am.\ Math.\ Soc., Providence, 1990.

\item{18.} Witten, E., {\it Nucl.\ Phys.}\ {\bf B322} 629 (1989);
{\bf B330}, 285 (1990).

\item{19.}  Martin, S.\ P., {\it Nucl.\ Phys.}\ {\bf B338}, 244 (1990).

\item{20.}  De Boer, J.\ and Goeree, J., {\it Commun.\ Math.\ Phys.}\
{\bf 139}, 267 (1991).

\item{21.} Gannon, T., {\it Nucl.\ Phys.}\ {\bf B396}, 708 (1993);

\item{}  Ruelle, Ph.,  Thiran, E., and  Weyers, J., {\it Nucl.\ Phys.}\
{\bf B402}, 693 (1993).

\item{22.} Guadagnini, E., {\it The Link Invariants of the Chern-Simons
Field Theory}, Walter De Gruyter, Berlin, 1993.

\item{23.} Chase, S.\ V.\ and Sweedler, M.\ E., {\it Hopf Algebras and
Galois Theory}, Lecture Notes in Math 97, Springer-Verlag, New York, 1969.

\item{24.} Grothendieck, A., {\it Esquisse d'un programme}, Rapport
Scientifique,  1984 (unpublished).

\item{25.} Drinfeld, V.\ G., {\it Lenin.\ Math.\ J.}\ {\bf 2}, 829 (1991).

\item{26.} Degiovanni, P., {\it Helv.\ Phys.\ Acta}\ {\bf 67}, 799 (1994).

\item{27.} Ihara, Y., in {\it Proc.\ of the ICM, Kyoto 1990}, Vol.\ I,
Springer-Verlag, Hong Kong, 1991.

\end